\documentclass[11pt]{article}

\usepackage{amsmath}
\usepackage{amsfonts}
\usepackage{url}
\usepackage{amsthm}
\usepackage{color}
\usepackage{stmaryrd}
\usepackage{graphicx}

\newtheorem{mylem}{Lemma}
\newtheorem{myass}{Assumption}

\newtheorem{myrem}{Remark}
\newtheorem{mythm}{Theorem}

\usepackage[margin=1.25cm]{geometry}
\begin{document}
\title{A priori error estimates of Adams--Bashforth  \\ discontinuous Galerkin methods for scalar nonlinear conservation laws}
\author{Charles Puelz and B\'eatrice Rivi\`ere}
\maketitle

\abstract{In this paper we show theoretical convergence of a second--order Adams--Bashforth discontinuous Galerkin method for approximating smooth solutions to scalar nonlinear conservation laws with E-fluxes.  A priori error estimates are also derived for a first--order forward Euler discontinuous Galerkin method.   Rates are optimal in time and suboptimal in space; they are valid under a CFL condition.%for both fully discrete schemes under a CFL condition, with the finite element space consisting of piecewise discontinuous polynomials of degree $k \geq 1$ for the forward Euler scheme and $k \geq 2$ for the Adams--Bashforth scheme.

\section{Introduction}

We consider approximating smooth solutions to the following nonlinear partial differential equation posed with initial conditions:
\begin{align}
\label{eq:PDE1}
&\frac{\partial u}{\partial t} + \frac{\partial}{\partial x}f(u) = s(u), \quad \text{in } \mathbb{R} \times (0,T],\\
\label{eq:PDE2}
&u = u_0, \quad \text{in } \in \mathbb{R}\times\{0\},
\end{align}
where $u:\mathbb{R} \times [0,T] \rightarrow \mathbb{R}$ and $f,s:\mathbb{R} \rightarrow \mathbb{R}$. The function $s$ is assumed to be Lipschitz. As typical for the numerical analysis of such problems \cite{ZS04, Zakerzadeh16}, we do not consider boundary conditions, and instead assume the solution has compact support in some interval $[0,L]$.  
%A discontinuous Galerkin scheme is employed for the spatial discretization, and we consider both a forward Euler and a second order Adams--Bashforth scheme for the time discretization.  

The focus of this work is the analysis of  the second order Adams--Bashforth method in time combined with the discontinuous Galerkin  method in space.  The main motivation for studying this discretization is its popularity in the hemodynamic modeling community for approximating a nonlinear hyperbolic system describing blood flow in an elastic vessel \cite{CK03}.  For a selection of work simulating this model with a discontinuous Galerkin spatial discretization coupled to the second order Adams--Bashforth scheme, see  \cite{SFPF03, SFPP03, Alastruey08, Boi15, WFL14, Dumas2017, Matthys2007, Bollache2014, Cascaval2016, PCRR16}.  To the best of our knowledge, there is little analysis for this fully discrete scheme.  The results presented in this paper for scalar hyperbolic equations provide a first step towards theoretically understanding the numerical approximation of the hyperbolic system modeling blood flow.
In addition, we provide an error analysis for the first order forward Euler in time combined with discontinuous Galerkin in space.

Discontinuous Galerkin schemes for hyperbolic conservations laws have been extensively studied, especially when coupled with Runge--Kutta methods for the time discretization.  This class of schemes was introduced in the series of papers by Cockburn, Shu, and co-authors \cite{CS_1, CS_2, CS_3, CS_4, CS_5}.  We recall the work from Zhang, Shu, and others analyzing Runge--Kutta discontinuous Galerkin methods applied to scalar conservation laws and symmetrizable systems \cite{ZS04, ZS06, ZS10, LSZ15}.  These papers establish error estimates for smooth solutions for both second and third order Runge--Kutta schemes.  Their analysis requires the CFL condition $\Delta t = O(h^{4/3})$ for the second order Runge--Kutta scheme and piecewise polynomials of degree two and higher.  The CFL condition $\Delta t = O(h)$ may be used for the third order Runge-Kutta scheme for piecewise polynomials of degree one and higher and for the second order Runge--Kutta scheme with piecewise linear polynomials.

Recent stability and convergence results have been obtained for IMEX (implicit--explicit) multistep schemes applied to a nonlinear convection diffusion equation, i.e. (\ref{eq:PDE1})--(\ref{eq:PDE2}) augmented with a nonzero diffusion term  \cite{Wang2016}.  These schemes implicitly discretize the diffusion term and explicitly discretize the hyperbolic term.  It is not immediately clear how to adapt the analysis to the case of zero diffusion since the estimates depend on the reciprocal of the diffusion parameter.

A summary of the paper is as follows.  In Section~\ref{sec:pb}, we introduce the numerical schemes, properties of the numerical flux, and inequalities related to projections.  The main results are also stated.  Section~\ref{sec:ABerror} and ~\ref{sec:FEerror} contain the proofs of the convergence results. In Section~\ref{sec:numer} we provide some numerical results for inviscid Burger's equation and a nonlinear hyperbolic system modeling blood flow in an elastic vessel.  Conclusions follow.

\section{Scheme and main results}
\label{sec:pb}

We define notation relevant for the spatial discretization of (\ref{eq:PDE1})--(\ref{eq:PDE2}) by the discontinuous Galerkin method.  To do this, we make a similar technical modification to the flux function as in \cite{ZS04}.  If the initial condition $u_0$ takes values within some open set $\Omega$, then locally in time the solution to (\ref{eq:PDE1})--(\ref{eq:PDE2}) also takes values in $\Omega$ \cite{Dafermos2010}.  We assume the flux function $f \in C^3(\mathbb{R})$ vanishes outside of $\Omega$ so derivatives up to third order are uniformly bounded, i.e. there exists some constant $C$ depending only on $f$ and its derivatives satisfying:
\begin{align}
\label{eq:fbounded}
	|f^{(\gamma)}(v)| \leq C, \quad \forall v \in \mathbb{R}, \quad \gamma = 1,2,3.
\end{align} 

Let the collection of intervals $\left( I_j \right)_{j=0}^{N}$ be a uniform partition of the interval $[0,L]$, with $I_j = [x_j,x_{j+1}]$ of  size $h$. Let $\mathbb{P}^k(I_j)$ denote the space of polynomials of degree $k$ on the interval $I_j$.  The approximation space is
\begin{equation}
\mathbb{V}_h=\{ \phi_h:[0,L] \rightarrow \mathbb{R} \text{ s.t. } \phi_h|_{I_j} \in \mathbb{P}^k(I_j),\quad\forall j = 0, \ldots, N \}.
\end{equation}
The space $L^2(0,L)$ is the standard $L^2$ space; let $(\cdot,\cdot)$ denote the $L^2$ inner-product over $\Omega$, with associated norm $\Vert \cdot\Vert$. 
Let $\Pi_h$ be the $L^2$ projection into $\mathbb{V}_h$:
\begin{equation}\label{eq:L2proj}
(\Pi_h v, \phi_h) = (v,\phi_h),\quad\forall\phi_h\in\mathbb{V}_h, \quad\forall v\in L^2(0,L).
\end{equation}
Define the notation for traces of a function $\phi:[0,L]\rightarrow \mathbb{R}$ to the boundaries of the intervals:
\begin{align}
\phi^\pm|_{x_j} &= \lim_{\varepsilon \rightarrow 0, \,  \varepsilon > 0} \phi(x_j \pm \varepsilon), \quad 1\leq j\leq N, \\
\phi^+|_{x_0} &= \lim_{\varepsilon \rightarrow 0, \, \varepsilon > 0} \phi(x_0 + \varepsilon), \\
\phi^-|_{x_{N+1}} &= \lim_{\varepsilon \rightarrow 0, \, \varepsilon > 0} \phi(x_{N+1} -\varepsilon).
\end{align}
The standard notation for jumps and averages at the interior nodes is given as follows:
\begin{align}
	[\phi]|_{x_j} &= \phi^-|_{x_j} - \phi^+|_{x_j}, \quad 1\leq j\leq N, \\
	\{\phi\}|_{x_j} &= \frac{1}{2}(\phi^-|_{x_j} + \phi^+|_{x_j}), \quad 1\leq j\leq N.
\end{align}
Let $\hat{f}$ denote the numerical flux, that is assumed to be Lipschitz continuous and consistent.
\begin{myass}
There is a constant $C_L > 0$ such that for any $p$, $q$, $u$, $v$ $\in \mathbb{R}$:
\begin{align}
        |\hat{f}(p, q) - \hat{f}(u, v)| \leq C_L \left( |p - u| + |q - v| \right),
\end{align}
and 
\begin{align}
\hat{f}(v,v) = f(v), \quad \forall v \in \mathbb{R}.
\end{align}
\label{myass1}
\end{myass}
We also assume that $\hat{f}$ belongs to the class of E--fluxes \cite{Di2011}. 
\begin{myass}
The numerical flux $\hat{f}$ is an E--flux, which means it satisfies, for all $w$ between $v^-$ and $v^+$, 
	\begin{equation}
		\left( \hat{f}(v^-,v^+) - f(w) \right) [v]|_{x_j} \geq 0, \quad 1 \leq j \leq N.
	\end{equation}
\label{myass2}
\end{myass}
An example of a numerical flux that satisfies Assumption~\ref{myass1} and Assumption~\ref{myass2} is the local Lax-Friedrichs flux, $\hat{f}_{LF}$,  defined by:
	\begin{align}
		\hat{f}_{LF}(v^-, v^+)|_{x_j} = \{ f(v) \}|_{x_j} + \frac{1}{2}J(v^-, v^+) [v]|_{x_j}, \quad \forall 1\leq j \leq N,
\label{eq:Lax}
	\end{align}
with 	
\begin{align}
J(v^-, v^+)|_{x_j} = \max_{\min(v^-|_{x_j}, v^+|_{x_j}) \leq w \leq \max(v^-|_{x_j}, v^+|_{x_j})} |f'(w)| , \quad \forall j = 1, \ldots, N.
\end{align}
Finally, we define a discrete function $\alpha$ at each interior node. The fact that $\alpha$ is nonnegative and
uniformly bounded is a key ingredient in the error analysis.
\begin{equation}
	\label{eq:alpha}
		\alpha(v)|_{x_j} =
		\left\{
		\begin{array}{lr}
		[v]^{-1}\left( \hat{f}(v^-, v^+) - f(\{v\}) \right)|_{x_j}, & \text{if } [v]|_{x_j} \neq 0, \\
		\frac{1}{2}\left|f'(\{v\}|_{x_j})\right|, & \text{if } [v]|_{x_j} = 0.
		\end{array} \right. 
	\end{equation}
\begin{mylem}
\label{lem:numflux}
There exist constants $C_\alpha$, $C_0$ and $C_1$ such that
\begin{align}
\label{eq:alpha0}
0 \leq \alpha(v)|_{x_j} &\leq C_\alpha, \quad  \forall (v^-, v^+) \in \mathbb{R}^2, \quad \forall \, 1\leq j\leq N,\\
	\label{eq:alpha1}
		\frac{1}{2} \Big|f'(\{v\}|_{x_j})\Big| &\leq \alpha(v)|_{x_j} + C_0 \Big|[v]|_{x_j}\Big|, \quad  \forall (v^-, v^+) \in \mathbb{R}^2, \quad \forall \,1\leq j\leq N, \\ 
	\label{eq:alpha2}
		\frac{1}{8} f''(\{v\}|_{x_j})[v]|_{x_j} &\leq \alpha(v)|_{x_j} + C_1 [v]^2|_{x_j}, \quad  \forall (v^-, v^+) \in \mathbb{R}^2, \quad \forall \,1\leq j\leq N.
	\end{align}
%[BR: we do not use \eqref{eq:alpha2}????]
The constants $C_0$ and $C_1$ depend on the derivatives of $f$.
\end{mylem}
The proof of Lemma~\ref{lem:numflux} follows the one in \cite{ZS04}; the definition for $\alpha$ slightly differs 
from the one given in \cite{ZS04} so that it is suitable for the error analysis of the
Adams--Bashforth scheme.

An additional assumption is made for the numerical flux. 
\begin{myass}
There is a constant $C > 0$ such that for any $v_h\in\mathbb{V}_h$ and $v\in \mathcal{C}(0,L)$:
\begin{align}
\vert \alpha(v_h)|_{x_j} -\alpha(v)|_{x_j} \vert \leq C \Vert v_h-v\Vert_\infty,\quad \forall 1\leq j\leq N.
\label{eq:ass3}
\end{align}
\label{myass3}
\end{myass}
\begin{myrem}
Assumption~\ref{myass3} is used in the error analysis for the Adams--Bashforth scheme.  It is easy to check
that the local Lax-Friedrichs flux defined by \eqref{eq:Lax} satisfies \eqref{eq:ass3}.
\end{myrem}

We now introduce the discontinuous Galerkin discretization on each interval.
\begin{align}
\mathcal{H}_j(v, \phi_h) &= \int_{I_j}f(v) \frac{d \phi_h}{dx} 
+ \int_{I_j}s(v) \phi_h -\hat{f}(v^-, v^+)|_{x_{j+1}} \phi_h^-|_{x_{j+1}} 
+ \hat{f}(v^-, v^+)|_{x_j} \phi_h^+|_{x_j} \quad \forall\, 1\leq j\leq N-1, \\
\mathcal{H}_0(v, \phi_h) &= \int_{I_0}f(v) \frac{d \phi_h}{dx} 
+ \int_{I_0}s(v) \phi_h -\hat{f}(v^-, v^+)|_{x_{1}} \phi_h^-|_{x_{1}}, \\
\mathcal{H}_{N}(v, \phi_h) &=\int_{I_N}f(v) \frac{d \phi_h}{dx} 
+ \int_{I_N}s(v) \phi_h + \hat{f}(v^-, v^+)|_{x_N} \phi_h^+|_{x_N}.
\end{align}
For some number $M>0$, define $\Delta t = T / M$. 
%We consider two schemes for the temporal discretization in this paper.  
The second order in time Adams--Bashforth scheme is:
given $u_h^0 \in\mathbb{V}_h$ and $u_h^1 \in\mathbb{V}_h$,  for $n = 1, \ldots, M-1$, seek $u_h^{n+1} \in \mathbb{V}_h$ satisfying
\begin{equation}
\label{eq:ABscheme}
 \int_{I_j}u_h^{n+1} \phi_h  =   \int_{I_j} u_h^{n} \phi_h  + \Delta t \frac{3}{2}\mathcal{H}_j(u_h^{n}, \phi_h) - \Delta t \frac{1}{2} \mathcal{H}_j(u_h^{n-1}, \phi_h),
\quad\forall \phi_h \in \mathbb{V}_h, \quad \forall\, 0\leq j\leq N.
\end{equation}
Since \eqref{eq:ABscheme} is a multi-step method, two starting values are needed.  We choose 
$u_h^0 = \Pi_h u_0$  for the initial value, and we choose $u_h^1 = \tilde{u}_h^1$ where
$\tilde{u}_h^1$  satisfies the first-order in time
forward Euler scheme defined below.

\noindent
%The standard first order in time forward Euler scheme is defined as:
With the choice $\tilde{u}_h^0 = \Pi_h u_0$, for $n = 0, \ldots, M-1$, seek $\tilde{u}_h^{n+1} \in \mathbb{V}_h$ satisfying
\begin{equation}
\label{eq:FEscheme}
 \int_{I_j}\tilde{u}_h^{n+1} \phi_h  =   \int_{I_j}\tilde{u}_h^{n} \phi_h  + \Delta t  \, \mathcal{H}_j(\tilde{u}_h^{n}, \phi_h),\quad
\forall \phi_h \in \mathbb{V}_h, \quad \forall\, 0\leq j\leq N.
\end{equation}
The initial value $u_h^1$ is computed using \eqref{eq:FEscheme} with a time step that is small enough so that the following
assumption holds:
\begin{equation}
\label{eq:uh1}
\| u_h^1 - \Pi_h u^1\| \leq h^{k+1/2}.
\end{equation}
Theorem~\ref{thm:FE} below shows that \eqref{eq:uh1} is a reasonable assumption if the time step used for
the forward Euler method is small enough.

The main result of this paper is the convergence result for the Adams-Bashforth scheme \eqref{eq:ABscheme}.
\begin{mythm}
\label{thm:AB}
Assume the exact solution $u$ belongs to $\mathcal{C}^2([0,T];H^{k+1}(\Omega))$. 
Let $u_h^1$ satisfy \eqref{eq:uh1}. 
Under Assumptions~\ref{myass1}, \ref{myass2}, \ref{myass3} and the CFL condition 
$\Delta t = O(h^2)$, there is a constant $C$ independent of $h$ and $\Delta t$ such that, 
 for $h$ sufficiently small, and for $k \geq 2$: 
	\begin{align}
		\max_{n = 0, \ldots, M}\|u^{n} - u_h^{n}\| \leq C (\Delta t^2 + h^{k+1/2}).
	\end{align}
\end{mythm}
The proof of Theorem~\ref{thm:AB} is given in Section~\ref{sec:ABerror}. An easy modification of the proof yields the
following convergence result for the forward Euler scheme \eqref{eq:FEscheme}. Its proof is outlined in Section~\ref{sec:FEerror}.
\begin{mythm}
\label{thm:FE}
Assume the exact solution $u$ belongs to $\mathcal{C}^2([0,T];H^{k+1}(\Omega))$. 
Let $(\tilde{u}_h^n)_n$ satisfy \eqref{eq:FEscheme}.
Under Assumptions~\ref{myass1}, \ref{myass2} and the CFL condition 
$\Delta t = O(h^2)$, for $h$ sufficiently small, and 
for $k \geq 1$, there is a constant $C$ independent of $h$ and $\Delta t$ such that:	
\begin{align}
	\max_{n = 0, \ldots, M}\|u^{n} - \tilde{u}_h^{n}\| \leq C (\Delta t + h^{k+1/2}).
\end{align}
\end{mythm}
\begin{myrem}
We remark that von Neumann stability analysis conducted in \cite{Deriaz2012} suggests a less restrictive CFL condition $\Delta t = O(h^{4/3})$ for the second order Adams--Bashforth scheme.  Our theoretical estimates require $\Delta t = O(h^2)$; at the moment we are unable to relax this condition. 
% Beatrice: what do you think?  I looked at the paper again, and they do some sort of stability analysis for scalar nonlinear conservation laws too.  maybe this can be useful.
\end{myrem}
 
\noindent
We finish this section by recalling inverse inequalities, trace inequalities and approximations results.
Let $\|v\|_\infty = \max_{x \in [0,L]}|v(x)|$ denote the sup-norm. There exists a constant $C$ independent of $h$ such that
\begin{align}
	\label{eq:inv1}
	\|\phi_h\|_\infty &\leq C h^{-1/2} \|\phi_h\|, \quad\forall \phi_h\in\mathbb{V}_h,\\
|\phi_h^{n,\pm}|_{x_j}| &\leq C h^{-1/2}\|\phi_h\|_{L^2(I_j)},\quad \forall \, 1\leq j\leq N, \quad\forall \phi_h\in\mathbb{V}_h,
	\label{eq:inv2}
\\
	\label{eq:inv3}	
	\left( \sum_{j = 0}^N\|\frac{d}{dx}\phi_h\|_{L^2(I_j)}^2 \right)^{1/2} &\leq C h^{-1} \|\phi_h\|, \quad\forall \phi_h\in\mathbb{V}_h.
\end{align}
For simplicity we denote $u^n$ the function $u$ evaluated at the time $t^n = n \Delta t$. The approximation error is denoted
\[
\eta^n  = u^n - \Pi_h u^n,
\]
and it satisfies the optimal a priori bounds
\begin{align}
\label{eq:est1}
	\|\eta^n\| &\leq C h^{k+1}, \\
|\eta^{n,\pm}|_{x_j}| & \leq C h^{k+1/2},\quad \forall \, 1\leq j\leq N,\\
%\label{eq:est2}
%	\|\eta^n\|_{\Gamma_h} &\leq C h^{k+1/2}, \\
\label{eq:est3}
	\|\eta^n\|_\infty &\leq C h^{k+1/2}, \\
\label{eq:est4}
	\|\eta^{n+1} - \eta^n \| &\leq C \Delta t\,h^{k+1}.
%\label{eq:est5}
	%\|\eta^{n+1} - \eta^n \|_{\Gamma_h} &\leq C \Delta t\,h^{k+1/2}. 
\end{align}
The constant $C$ is independent of $h, \Delta t$ but depends on the exact solution $u$ and its derivatives.

%%HERE
\section{Proof of Theorem~\ref{thm:AB}}
\label{sec:ABerror}

For the error analysis, we denote
\[
\chi^n = u_h^n-\Pi_h u^n.
\]
The proof of Theorem~\ref{thm:AB} is based on an induction hypothesis:
\begin{align}
	\|\chi^\ell\| \leq h^{3/2}, \quad \forall  0\leq \ell \leq M.
\label{eq:inductAB1}
\end{align}
Since $\chi^0=0$, the hypothesis (\ref{eq:inductAB1}) is trivially satisfied for $\ell=0$.  
With the assumption \eqref{eq:uh1}, it is also true for $\ell=1$.
Fix $\ell \in \{2,\ldots,M\}$ and assume that 
\begin{align}
	\|\chi^n\| \leq h^{3/2}, \quad \forall  0\leq n \leq \ell-1.
\label{eq:inductAB2}
\end{align}
We will show that (\ref{eq:inductAB2}) is valid for $n = \ell $.  We begin by deriving an error inequality.
We fix an interval $I_j$ for $0\leq j\leq N$. 
It is easy to see that the scheme is consistent in space and the exact solution satisfies
\begin{align}
	\frac{3}{2}\int_{I_j}& u_t^n \phi_h - \frac{1}{2}\int_{I_j}u_t^{n-1} \phi_h = \frac{3}{2}\mathcal{H}_j(u^n, \phi_h) - \frac{1}{2}\mathcal{H}_j(u^{n-1}, \phi_h), \quad \forall 1\leq n \leq M-1.
\label{eq:exactcons}
\end{align}
In the above, the notation $u_t^n$ is used for the time derivative of $u$ evaluated at $t^n$.
Subtracting \eqref{eq:exactcons} from (\ref{eq:ABscheme}) and rearranging terms, one obtains:
\begin{align*}
	\int_{I_j}&\left(u_h^{n+1} -u_h^{n} - \Delta t\frac{3}{2} u_t^n  + \Delta t\frac{1}{2} u_t^{n-1} \right)\phi_h \nonumber \\
	 &= \Delta t\frac{3}{2}(\mathcal{H}_j(u_h^{n}, \phi_h) - \mathcal{H}_j(u^{n}, \phi_h) ) - \Delta t\frac{1}{2} (\mathcal{H}_j(u_h^{n-1}, \phi_h) - \mathcal{H}_j(u^{n-1}, \phi_h) ), \quad \forall 1\leq n \leq M-1.
\end{align*}
Summing over the elements $j = 0, \ldots, N$ and adding and subtracting the $L^2$ projection 
of $u$ at $t^n$ and $t^{n+1}$ yields the equality:
\begin{align}
\label{eq:foo54}
	\int_{0}^L (\chi^{n+1} - \chi^n) \phi_h &=  \int_{0}^L \left(u^{n} -u^{n+1} + \Delta t \frac{3}{2} u_t^n  - \Delta t\frac{1}{2} u_t^{n-1}  \right)\phi_h  +  \int_{0}^L (\eta^{n+1} - \eta^n) \phi_h  + b^n(\phi_h),
\end{align}
with the following definition for $n \geq 1$ 
\begin{align*}
b^n(\phi_h) &= \Delta t \frac{3}{2} \sum_{j = 0}^N \left( \mathcal{H}_j(u_h^{n}, \phi_h) - \mathcal{H}_j(u^{n}, \phi_h)\right) - \Delta t  \frac{1}{2} \sum_{j = 0}^N \left( \mathcal{H}_j(u_h^{n-1}, \phi_h) - \mathcal{H}_j(u^{n-1}, \phi_h) \right).
\end{align*}
The second term on the right hand side of (\ref{eq:foo54}) vanishes due to the property \eqref{eq:L2proj} of the local $L^2$ projection.  To handle the first term, we obtain from the following Taylor expansions for some $\tilde{\zeta} \in [t^{n-1},t^n]$ and some $\zeta \in [t^{n}, t^{n+1}]$:
\begin{align*}
	u^{n+1} - u^n &= \Delta t u_t^n  + \frac{1}{2} \Delta t^2 u_{tt}^n  + \frac{1}{6}\Delta t^3 u_{ttt}|_\zeta, \\
	u_t^{n-1} - u_t^n  &= -\Delta t u_{tt}^n  + \frac{1}{2} \Delta t^2 u_{ttt}|_{\tilde{\zeta}}.
\end{align*}
Thus we have
\begin{align*}
	u^{n} -u^{n+1} + \Delta t \frac{3}{2} u_t^n  - \Delta t\frac{1}{2} u_t^{n-1} &= - \Delta t^3 (\frac{1}{6}u_{ttt}|_\zeta +\frac{1}{4}u_{ttt}|_{\tilde{\zeta}}).
\end{align*}
Hence (\ref{eq:foo54}) becomes:
\begin{align}
\label{eq:foo10}
 	\int_{0}^L \left(\chi^{n+1} - \chi^n \right) \phi_h  \leq C \Delta t^3 \int_{0}^L |\phi_h| +  b^n(\phi_h).
 \end{align}
Cauchy Schwarz's inequality and Young's inequalities imply:
\begin{align}
\label{eq:ABerrorequation}
 	\int_{0}^L \left(\chi^{n+1} - \chi^n \right) \phi_h  \leq C \Delta t^5 +  \Delta t  \|\phi_h\|^2   + b^n(\phi_h).
\end{align}
We choose $\phi_h = \chi^n$ in inequality (\ref{eq:ABerrorequation}) to obtain:
\begin{align*}
 	\int_{0}^L \left(\chi^{n+1} - \chi^n \right) \chi^n  \leq  C \Delta t^5 +  \Delta t  \|\chi^n\|^2  + b^n(\chi^n).	
\end{align*} 
So, the following error inequality holds for $n \geq 1$:
\begin{align}
\label{eq:ABerror}
	\frac{1}{2} \| \chi^{n+1} \|^2 - \frac{1}{2} \| \chi^{n} \|^2  \leq  
 C \Delta t^5 +  \Delta t  \|\chi^n\|^2  
+\frac{1}{2}\| \chi^{n+1} - \chi^{n} \|^2 
+ b^n(\chi^n).
\end{align}
It remains to handle the last two terms in \eqref{eq:ABerror}.  The proofs of the following two lemma are given in the next section.
\begin{mylem}
\label{cor:AB1} 
Assume that $\Delta t = O(h^2)$.  The following holds for $n \geq 1$:
\begin{align}
\label{eq:cor:AB1}
	\|\chi^{n+1} - \chi^n\|^2 &\leq C \Delta t^6 + C \Delta t \, ( \|\chi^n\|^2 + \|\chi^{n-1}\|^2) + C \Delta t \, h^{2k+2}.
\end{align}
\end{mylem}
\begin{mylem}
\label{cor:AB2} 
Let $n\geq 2$ and assume $\|\chi^n\| \leq h^{3/2}$, $\Vert \chi^{n-1}\Vert \leq h^{3/2}$,  and $\Delta t = O(h^2)$.  The following holds:
\begin{align}
\label{eq:cor:AB2}
b^n(\chi^n) \leq &C \Delta t (\Vert \chi^{n}\Vert^2+\Vert\chi^{n-1}\Vert^2) 
+ C \Delta t^6 
+ C \Delta t\, (1+ 2 \varepsilon^{-1}) h^{2k+1}\nonumber\\
&-(\frac12-2 \varepsilon) \Delta t \sum_{j=1}^N \alpha(u_h^n)|_{x_j} [\chi^n]^2|_{x_j}
-(\frac12-2 \varepsilon) \Delta t \sum_{j=1}^N \alpha(u_h^{n-1})|_{x_j} [\chi^{n-1}]^2|_{x_j}, \quad
\forall\varepsilon >0.
\end{align}
For $n = 1$ one has the following:
\begin{align}
\label{eq:cor:AB2n1}
b^1(\chi^1) \leq &C \Delta t \Vert \chi^{1}\Vert^2
+ C \Delta t\, (1+ 2 \varepsilon^{-1}) h^{2k+1}\nonumber\\
&+3\Vert \chi^1\Vert^2-(\frac12-2 \varepsilon) \Delta t \sum_{j=1}^N \alpha(u_h^1)|_{x_j} [\chi^1]^2|_{x_j},\quad
\forall\varepsilon >0.
\end{align}
\end{mylem}
Substituting the bounds from \eqref{eq:cor:AB1}, \eqref{eq:cor:AB2}, \eqref{eq:cor:AB2n1} (with $\varepsilon = 1/4$), and using the fact that 
$\alpha(u_h^n)$ and $\alpha(u_h^{n-1})$ are nonnegative, the error inequality \eqref{eq:ABerror} simplifies to:
\begin{align}
\| \chi^{n+1} \|^2 - \| \chi^{n} \|^2 &\leq C \Delta t^5 +  C \Delta t \, ( \|\chi^n\|^2 + \|\chi^{n-1}\|^2 + \Vert \chi^{n-2}\Vert^2)  + C \Delta t\,h^{2k+1}, \quad
n\geq 2,
\label{eq:intern}
\end{align}
and
\begin{align}
\| \chi^{n+1} \|^2 - \| \chi^{n} \|^2 &\leq C \Delta t^5 +  C \Delta t \,  \|\chi^n\|^2   + C \Delta t\,h^{2k+1}
+ C \Vert \chi^n\Vert^2, \quad
n = 1.
\label{eq:intern1}
\end{align}
Summing \eqref{eq:intern} from $n = 2, \ldots, \ell-1$ and adding to \eqref{eq:intern1} one obtains:
\begin{align*}
  \| \chi^{\ell} \|^2  &\leq C \Delta t^4 + Ch^{2k+1} +  C\|\chi^1\|^2  +  C \Delta t \, \sum_{n = 0}^{\ell-1} \|\chi^n\|^2. 
  \end{align*}
Gronwall's inequality  and assumption \eqref{eq:uh1} immediately gives
\begin{align*}
	\| \chi^\ell \|^2 \leq C_2 T \mathrm{e}^T  \left( \Delta t^4 +  h^{2k+1} \right),
\end{align*}
where $C_2$ is independent of $\ell$, $h$ and $\Delta t$.  Employing the CFL condition $\Delta t = O(h^2)$, one has:
\begin{align*}
	\| \chi^\ell\| \leq \left(C_2 T \mathrm{e}^T \right)^{1/2} \left( h^4 +  h^{k+1/2} \right).
\end{align*}
The induction proof is complete if $h$ is small enough so that
\begin{align*}
	C_2 T \mathrm{e}^T h < 1,
\end{align*} 
implying that for $k \geq 2$:
\begin{align*}
	\| \chi^\ell \| 
	\leq  \left(C_2 T \mathrm{e}^T \right)^{1/2} h \left( h^3 +  h^{k-1/2} \right)  \leq h^{3/2}.  
\end{align*}
Since $\|\eta^n\| \leq C h^{k+1}$ and $\|u^n - u_h^n\| \leq \|\eta^n\| + \|\chi^n\|$ we can conclude:
\begin{align*}
	\|u^{n} - u_h^{n}\| \leq C\left(\Delta t^2 + h^{k+1/2} \right).
\end{align*}

\subsection{Proof of Lemma~\ref{cor:AB1}}

Choose $\phi_h = \chi^{n+1} - \chi^n$ in (\ref{eq:foo10}) and use Cauchy-Schwarz's and Young's inequalities to obtain:
\begin{align}
	\|\chi^{n+1} - \chi^n\|^2 \leq  C\Delta t^6 + 2\, b^n(\chi^{n+1} - \chi^n).
\label{eq:chinn1}
\end{align}	
We will now obtain a bound for $b(\phi_h)$ for any $\phi_h \in \mathbb{V}_h$. 
By definition, we write
\[
b^n(\phi_h) = \theta_1 + \theta_2 + \theta_3,
\]
where
\begin{eqnarray}
\theta_1  = \frac32 \Delta t \sum_{j=0}^N \int_{I_j} (f(u_h^n)-f(u^n))\frac{d\phi_h}{dx}
-\frac12 \Delta t \sum_{j=0}^N \int_{I_j} (f(u_h^{n-1})-f(u^{n-1}))\frac{d\phi_h}{dx}\nonumber\\
-\frac32 \Delta t \sum_{j=1}^N (f(\{u_h^n\})-f(u^n))|_{x_j} [\phi_h]|_{x_j}
+\frac12 \Delta t \sum_{j=1}^N (f(\{u_h^{n-1}\})-f(u^{n-1}))|_{x_j} [\phi_h]|_{x_j},
\label{eq:theta1}
\end{eqnarray}
\begin{eqnarray}
\label{eq:theta2}
\theta_2 =  \Delta t \sum_{j=0}^N \int_{I_j}  \left(\frac32 (s(u_h^n)-s(u^n))-\frac12 (s(u_h^{n-1})-s(u^{n-1}))\right)\phi_h,
\end{eqnarray}
\begin{eqnarray}
\label{eq:theta3}
\theta_3 = -\frac32 \Delta t \sum_{j=1}^N (\hat{f}(u_h^{n,-},u_h^{n,+})-f(\{u_h^n\}))|_{x_j} [\phi_h]|_{x_j}
+\frac12 \Delta t \sum_{j=1}^N (\hat{f}(u_h^{n-1,-},u_h^{n-1,+})-f(\{u_h^{n-1}\}))|_{x_j} [\phi_h]|_{x_j}.
\end{eqnarray}
Using Taylor expansions, we write for some $\zeta_1^n$, $\zeta_2^n$, $\zeta_1^{n-1}$ and $\zeta_2^{n-1}$:
\begin{align*}
	f(u_h^n) - f(u^n) &= f'(\zeta_1^n)(u_h^n - u^n) = f'(\zeta_1^n)(\chi^n-\eta^n), \nonumber\\
	f(\{u_h^n\}) - f(u^n) &= f'(\zeta_2^n)(\{u_h^n\} - \{u^n\}) =  f'(\zeta_2^n)(\{\chi^n\} - \{\eta^n\}), \nonumber\\
	f(u_h^{n-1}) - f(u^{n-1}) &= f'(\zeta_1^{n-1})(u_h^{n-1} - u^{n-1}) = f'(\zeta_1^{n-1})(\chi^{n-1}-\eta^{n-1}), \nonumber\\
	f(\{u_h^{n-1}\}) - f(u^{n-1}) &= f'(\zeta_2^n)(\{u_h^{n-1}\} - \{u^{n-1}\}) =  f'(\zeta_2^{n-1})(\{\chi^{n-1}\} - \{\eta^{n-1}\}).
\end{align*}
Using the above expansions in the definition of $\theta_1$, trace inequalities and the CFL condition $\Delta t =\mathcal{O}(h^2)$, we
can obtain for any $\varepsilon > 0$
\begin{equation}
|\theta_1|  \leq \varepsilon\Vert \phi_h\Vert^2 + C \varepsilon^{-1}  \Delta t (\Vert \chi^n\Vert^2 +\Vert \chi^{n-1}\Vert^2)
+C \varepsilon^{-1} \, \Delta t \, h^{2k+2}.
\label{eq:lem1theta1}
\end{equation}
The term $\theta_2$ is bounded using Lipschitz continuity of $s$, approximation results, Cauchy-Schwarz's and Young's inequalities.
For any $\varepsilon>0$, we have
\begin{align*}
	\theta_2 
\leq C \varepsilon^{-1} \Delta t^2 h^{2k+2} + C\varepsilon^{-1} \Delta t^2 (\|\chi^n\|^2+\Vert \chi^{n-1}\Vert^2)  
+ \varepsilon \|\phi_h\|^2.
\end{align*}
Lastly, the term $\theta_3$ can be rewritten using the definition  \eqref{eq:alpha}. 
\begin{align*}
	\theta_3 = & -\frac32 \Delta t \sum_{j = 1}^N \alpha(u_h^n)|_{x_j} [u_h^n]|_{x_j} [\phi_h]|_{x_j} 
+\frac12 \Delta t \sum_{j = 1}^N \alpha(u_h^{n-1})|_{x_j} [u_h^{n-1}]|_{x_j} [\phi_h]|_{x_j}
\nonumber\\
= & -\frac32 \Delta t \sum_{j = 1}^N \alpha(u_h^n)|_{x_j} [\chi^n-\eta^n]|_{x_j} [\phi_h]|_{x_j} 
+\frac12 \Delta t \sum_{j = 1}^N \alpha(u_h^{n-1})|_{x_j} [\chi^{n-1}-\eta^{n-1}]|_{x_j} [\phi_h]|_{x_j}.
\end{align*}
Using Young's and Cauchy-Schwarz's inequalities, approximation results, trace inequalities, boundedness of $\alpha$ and the CFL condition, we have
\[
|\theta_3 | \leq \varepsilon \Vert \phi_h\Vert^2 
+ C \varepsilon^{-1} \Delta t (\Vert \chi^n\Vert^2+\Vert \chi^{n-1}\Vert^2)
+ C \varepsilon^{-1} \Delta t \, h^{2k+2}.
\]
Combining the bounds above yields
\begin{equation}
b(\phi_h) \leq \varepsilon \Vert \phi_h\Vert^2 +
C \varepsilon^{-1} \Delta t (\Vert \chi^n\Vert^2+\Vert \chi^{n-1}\Vert^2)
+ C \varepsilon^{-1} \Delta t \, h^{2k+2}, \quad \forall \varepsilon >0, \quad \forall\phi_h\in \mathbb{V}_h.
\label{eq:bphihbound}
\end{equation}
We choose $\varepsilon = 1/4$ and $\phi_h  = \chi^n-\chi^{n-1}$ in \eqref{eq:bphihbound} and substitute the bound 
in \eqref{eq:chinn1} to obtain \eqref{eq:cor:AB1}.
\begin{align}
	\|\chi^{n+1} - \chi^n\|^2 \leq  C\Delta t^6 
+ C  \Delta t (\Vert \chi^n\Vert^2+\Vert \chi^{n-1}\Vert^2)
+ C \Delta t \, h^{2k+2}. 
\label{eq:chinn2}
\end{align}

\subsection{Proof of Lemma~\ref{cor:AB2}}

As in the proof of Lemma~\ref{cor:AB1}, we write 
\[
b^n(\chi^n) = \theta_1 + \theta_2 + \theta_3,
\]
where the definitions of $\theta_1, \theta_2, \theta_3$ are given in \eqref{eq:theta1}, \eqref{eq:theta2} and \eqref{eq:theta3} respectively
for the particular choice $\phi_h = \chi^n$.  Unfortunately we cannot make use of the bound \eqref{eq:bphihbound} since the factor
$\Delta t$ is missing in front of $\varepsilon \Vert \phi_h \Vert^2$.  A more careful analysis is needed, and we will take
advantage of the CFL condition.
%%%%
Define
\begin{equation}
\mathcal{F}(n,\phi_h) = 
\Delta t \sum_{j=0}^N \int_{I_j} (f(u_h^n)-f(u^n))\frac{d\phi_h}{dx}
- \Delta t \sum_{j=1}^N (f(\{u_h^n\})-f(u^n))|_{x_j} [\phi_h]|_{x_j}.
\label{eq:defFn}
\end{equation}
Using the function $\mathcal{F}$ which is linear in its second argument, we rewrite the term $\theta_1$ as
\[
\theta_1 = \frac32 \mathcal{F}(n,\chi^n) -\frac12\mathcal{F}(n-1,\chi^{n-1}) + \frac12\mathcal{F}(n-1,\chi^{n-1}-\chi^n).
\]
We now state a bound for the term $\mathcal{F}(n,\chi^n)$. 
\begin{equation}
\mathcal{F}(n,\chi^n)\leq C \Delta t \Vert \chi^n\Vert^2
+ C(1+\varepsilon^{-1}) \Delta t \, h^{2k+1}
+\varepsilon \Delta t \sum_{j=1}^N \alpha(u_h^n)|_{x_j} [\chi^n]^2|_{x_j},\quad\forall \varepsilon > 0.
\label{eq:boundF}
\end{equation}
The proof of \eqref{eq:boundF} is technical and can be found in Appendix~\ref{app:boundF}.
The bound for $\mathcal{F}(n-1,\chi^{n-1})$ is identical.
\begin{equation}
\mathcal{F}(n-1,\chi^{n-1})\leq C \Delta t \Vert \chi^{n-1}\Vert^2
+ C(1+\varepsilon^{-1}) \Delta t \, h^{2k+1}
+\varepsilon \Delta t \sum_{j=1}^N \alpha(u_h^{n-1})|_{x_j} [\chi^{n-1}]^2|_{x_j},\quad\forall \varepsilon > 0.
\label{eq:boundFn}
\end{equation}
We are left with bounding $\mathcal{F}(n-1,\chi^{n-1}-\chi^n)$. Following the technique used for bound \eqref{eq:lem1theta1}, we can obtain
\begin{eqnarray}
\mathcal{F}(n-1,\chi^{n-1}-\chi^n)   \leq \Vert \chi^{n-1}-\chi^n\Vert^2 
+ C  \Delta t \Vert \chi^{n-1}\Vert^2 +C  \Delta t \, h^{2k+2}.
\label{eq:boundFnn}
\end{eqnarray}
Combining the above with \eqref{eq:chinn2}, we have for $n\geq 2$
\begin{align*}
\theta_1 \leq  & C \Delta t (\Vert \chi^n\Vert^2+ \Vert \chi^{n-1}\Vert^2 + \Vert \chi^{n-2}\Vert^2)
+ C (1+2 \varepsilon^{-1}) \Delta t \, h^{2k+1}
\nonumber\\
 & +\varepsilon \Delta t \sum_{j=1}^N \alpha(u_h^n)|_{x_j} [\chi^n]^2|_{x_j}
 +\varepsilon \Delta t \sum_{j=1}^N \alpha(u_h^{n-1})|_{x_j} [\chi^{n-1}]^2|_{x_j}
+ C \Delta t^6, \quad\forall \epsilon>0.
\end{align*}
For $n=1$, since $\chi^0 = 0$, inequalities \eqref{eq:boundF} and \eqref{eq:boundFnn} imply
\begin{align*}
\theta_1 \leq   C \Delta t \Vert \chi^1\Vert^2
+ C (1+\varepsilon^{-1}) \Delta t \, h^{2k+1}
  +\varepsilon \Delta t \sum_{j=1}^N \alpha(u_h^1)|_{x_j} [\chi^1]^2|_{x_j}
+ \Vert \chi^1\Vert^2,
\quad\forall \epsilon>0.
\end{align*}

%%%%
The term $\theta_2$ is bounded using Lipschitz continuity of $s$, approximation results, Cauchy-Schwarz's
inequality:
\[
\theta_2 \leq C \Delta t (\Vert \chi^n\Vert+\Vert \eta^n\Vert+\Vert\chi^{n-1}\Vert+\Vert\eta^{n-1}\Vert) \Vert\chi^n\Vert
\leq C \Delta t (\Vert\chi^n\Vert^2 +\Vert\chi^{n-1}\Vert^2) + C \Delta t \, h^{2k+2}.
\]
%%%
For the term $\theta_3$, we use the definition \eqref{eq:alpha} and write
\[
\theta_3 = -\frac32 \Delta t \sum_{j = 1}^N \alpha(u_h^n)|_{x_j} [\chi^n-\eta^n]|_{x_j} [\chi^n]|_{x_j} 
+\frac12 \Delta t \sum_{j = 1}^N \alpha(u_h^{n-1})|_{x_j} [\chi^{n-1}-\eta^{n-1}]|_{x_j} [\chi^n]|_{x_j}.
\]
After some manipulation we rewrite $\theta_3$ as:
\begin{eqnarray}
\theta_3 = -\frac12 \Delta t \sum_{j = 1}^N \alpha(u_h^n)|_{x_j} [\chi^n]^2|_{x_j}
-\frac12 \Delta t \sum_{j = 1}^N \alpha(u_h^{n-1})|_{x_j} [\chi^{n-1}]^2|_{x_j}\nonumber\\
+\Delta t \sum_{j=1}^N \left(\alpha(u_h^{n-1})-\alpha(u_h^n)\right)|_{x_j} [\chi^{n-1}-\eta^{n-1}]|_{x_j} [\chi^{n-1}]|_{x_j}
-\frac12\Delta t \sum_{j=1}^N  \alpha(u_h^{n-1})|_{x_j} [\chi^{n-1}-\eta^{n-1}]|_{x_j} [\chi^{n-1}-\chi^{n}]|_{x_j}
\nonumber\\
+\Delta t \sum_{j=1}^N  \alpha(u_h^n)|_{x_j} [\chi^{n-1}-\eta^{n-1}]|_{x_j} [\chi^{n-1}-\chi^{n}]|_{x_j}
+\Delta t \sum_{j=1}^N  \alpha(u_h^n)|_{x_j} [(\chi^{n-1}-\chi^{n})-(\eta^{n-1}-\eta^{n})]|_{x_j} [\chi^{n}]|_{x_j}
\nonumber \\
+\frac12 \Delta t \sum_{j=1}^N  \alpha(u_h^n)|_{x_j} [\eta^n]|_{x_j} [\chi^n]|_{x_j} 
+\frac12 \Delta t \sum_{j=1}^N  \alpha(u_h^{n-1})|_{x_j} [\eta^{n-1}]|_{x_j} [\chi^{n-1}]|_{x_j}. 
\label{eq:inttheta3}
\end{eqnarray}

We now bound the terms in the right-hand side of \eqref{eq:inttheta3} except for the first two terms.
We write
\[
\alpha(u_h^{n-1})|_{x_j} - \alpha(u_h^{n})|_{x_j}
=(\alpha(u_h^{n-1})|_{x_j} - \alpha(u^{n-1})|_{x_j})
+(\alpha(u^{n-1})|_{x_j} - \alpha(u^{n})|_{x_j})
-(\alpha(u_h^n)|_{x_j}-\alpha(u^n)|_{x_j}).
\]
From \eqref{eq:ass3} and \eqref{eq:alpha}, we have
\[
|\alpha(u_h^{n-1})|_{x_j} - \alpha(u_h^{n})|_{x_j}|
\leq C\Vert u_h^{n-1} - u^{n-1}\Vert_{\infty} + C \Vert u_h^n-u^n\Vert_{\infty}
+\frac12\left|  \, \vert f'(u^{n-1})|_{x_j}\vert -\vert f'(u^n)|_{x_j}\vert\, \right|.
\]
With a Taylor expansion, we obtain
\begin{align*}
	\Big| \alpha(u_h^{n-1})|_{x_j} &- \alpha(u_h^{n})|_{x_j} \Big| \leq  C \left(\|u^{n-1} - u_h^{n-1}\|_\infty + \|u^{n} - u_h^n\|_\infty + \Delta t \right), \quad \forall 1\leq j\leq N.
\end{align*}
With the assumption $\Vert \chi^n\Vert \leq h^{3/2}$ and $\Vert \chi^{n-1}\Vert \leq h^{3/2}$,  bound \eqref{eq:inv1}
and approximation results, we have 
\begin{align*}
	\Big| \alpha(u_h^{n-1})|_{x_j} &- \alpha(u_h^{n})|_{x_j} \Big| \leq  C (h + \Delta t), \quad \forall 1\leq j \leq N.
\end{align*}
Using trace inequalities, we then have
\[
\Delta t \sum_{j=1}^N \left(\alpha(u_h^{n-1})|_{x_j}-\alpha(u_h^n)|_{x_j}\right)[\chi^{n-1}]^2|_{x_j}
\leq C \Delta t (1 + h^{-1} \Delta t) \Vert \chi^{n-1}\Vert^2. 
\]
With the CFL condition, we conclude
\[
\Delta t \sum_{j=1}^N \left(\alpha(u_h^{n-1})|_{x_j}-\alpha(u_h^n)|_{x_j}\right)[\chi^{n-1}]^2|_{x_j}
\leq C \Delta t \Vert \chi^{n-1} \Vert^2.
\]
Similarly we have
\[
-\Delta t \sum_{j=1}^N \left(\alpha(u_h^{n-1})|_{x_j}-\alpha(u_h^n)|_{x_j}\right)[\eta^{n-1}]|_{x_j}[\chi^{n-1}]|_{x_j}
\leq C \Delta t \Vert \chi^{n-1} \Vert^2 + C \Delta t \, h^{2k+1}.
\]
The fourth term in \eqref{eq:inttheta3} is bounded by Cauchy-Schwarz's inequality,  trace inequalities, approximation results, the 
CFL condition and \eqref{eq:alpha0}:
\begin{align*}
\frac12 \Delta t \sum_{j=1}^N \alpha(u_h^{n-1})|_{x_j} [\chi^{n-1}-\eta^{n-1}]|_{x_j}[\chi^{n-1}-\chi^n]|_{x_j}
&\leq \Vert \chi^{n-1}-\chi^n\Vert^2 + C \Delta t^2 h^{-2} \Vert \chi^{n-1}\Vert^2 + C \Delta t^2 \, h^{2k} \nonumber\\
&\leq \Vert \chi^{n-1}-\chi^n\Vert^2  +  C \Delta t \Vert \chi^{n-1}\Vert^2 + C \Delta t h^{2k+2}.
\end{align*}
The fifth term in \eqref{eq:inttheta3} is handled exactly like the fourth term.
Similarly the first part in the sixth term has the following bound:
\[
\Delta t \sum_{j=1}^N  \alpha(u_h^n)|_{x_j} [\chi^{n-1}-\chi^{n}]|_{x_j} [\chi^{n}]|_{x_j}
\leq \Vert \chi^{n-1}-\chi^n\Vert^2 + C \Delta t \Vert \chi^n\Vert^2.
\]
For the second part, we use a Taylor expansion in time and the CFL condition:
\[
\Delta t \sum_{j=1}^N  \alpha(u_h^n)|_{x_j} [\eta^{n-1}-\eta^{n}]|_{x_j} [\chi^{n}]|_{x_j}
\leq C \Delta t^2 h^{k} \Vert \chi^n\Vert
\leq C \Delta t \Vert \chi^n\Vert^2 + C \Delta t \, h^{2k+2}.
\]
The last two terms in \eqref{eq:inttheta3} are treated almost identically, using approximation results, and
the boundedness of $\alpha$:
\begin{align*}
\frac12 \Delta t \sum_{j=1}^N \alpha(u_h^n)|_{x_j} [\eta^n]|_{x_j} [\chi^n]|_{x_j}
&+\frac12 \Delta t \sum_{j=1}^N \alpha(u_h^{n-1})|_{x_j} [\eta^{n-1}]|_{x_j} [\chi^{n-1}]|_{x_j}
\leq C \varepsilon^{-1} \Delta t \, h^{2k+1} 
+ \varepsilon \Delta t \sum_{j=1}^N \alpha(u_h^n)|_{x_j} [\chi^n]^2|_{x_j}\nonumber\\
&+ \varepsilon \Delta t \sum_{j=1}^N \alpha(u_h^{n-1})|_{x_j} [\chi^{n-1}]^2|_{x_j},\quad
\forall \varepsilon>0.
\end{align*}
To summarize, with \eqref{eq:chinn2}, the term $\theta_3$ is bounded as:
\begin{align*}
\theta_3 \leq &C \Delta t (\Vert \chi^{n}\Vert^2+\Vert\chi^{n-1}\Vert^2) 
+ C \Delta t^6 
+ C \Delta t\, (1+\varepsilon^{-1}) h^{2k+1}\nonumber\\
&-(\frac12-\varepsilon) \Delta t \sum_{j=1}^N \alpha(u_h^n)|_{x_j} [\chi^n]^2|_{x_j}
-(\frac12-\varepsilon) \Delta t \sum_{j=1}^N \alpha(u_h^{n-1})|_{x_j} [\chi^{n-1}]^2|_{x_j}, \quad
\forall\varepsilon >0, \quad n\geq 2.
\end{align*}
For $n=1$, the term $\theta_3$ is simply bounded as:
\begin{align*}
\theta_3 \leq C \Delta t \Vert \chi^{1}\Vert^2
+ C \Delta t\, (1+\varepsilon^{-1}) h^{2k+1}
-(\frac12-\varepsilon) \Delta t \sum_{j=1}^N \alpha(u_h^1)|_{x_j} [\chi^1]^2|_{x_j}
+ 2\Vert \chi^1\Vert^2,
\forall\varepsilon >0.
\end{align*}

Combining the bounds above for $\theta_i$, $1\leq i\leq 3$, we conclude the proof.

\section{Proof of Theorem~\ref{thm:FE}}
\label{sec:FEerror}

The proof for the forward Euler scheme is also done by induction.  It is a less technical proof than
for the Adams--Bashforth scheme.  We skip many details and give an outline of the proof.
Denote
\[
\xi^n = \tilde{u}_h^n-\Pi_h u^n.
\]
The induction hypothesis is less restrictive than for the Adams-Bashforth method, which yields
a convergence result that is valid for polynomials of degree one and above.
\begin{align}
	\|\xi^\ell\| \leq h, \quad \forall  0\leq \ell \leq M.
\label{eq:induct1}
\end{align}
Since $\xi^0=0$, the hypothesis (\ref{eq:induct1}) is trivially satisfied for $\ell=0$.  
Fix $\ell \in \{1,\ldots,M\}$ and assume that 
\begin{align}
	\|\xi^n\| \leq h, \quad \forall  0\leq n \leq \ell-1.
\label{eq:induct2}
\end{align}
We now have to  show that (\ref{eq:induct2}) is valid for $n = \ell $.  We begin by deriving an error inequality.
We fix an interval $I_j$ for $0\leq j\leq N$.  Using consistency in space of the scheme:
\begin{align}
\label{eq:exact_scheme}
	\int_{I_j}& u_t^n \phi_h = \mathcal{H}_j(u^n, \phi_h), \quad 0\leq n\leq M, 
\end{align}
we obtain, after some manipulation, the error equation:
\begin{align}
\label{eq:foo45}
	& \int_{I_j}\left(\xi^{n+1} - \xi^n\right) \phi_h = \int_{I_j}\left(\Delta t \, u_t^n  - u^{n+1} + u^n \right) \phi_h + \int_{I_j} \left(\eta^{n+1} - \eta^n \right)  \phi_h + \Delta t\left(\mathcal{H}_j(u_h^{n}, \phi_h) - \mathcal{H}_j(u^n, \phi_h)\right).
\end{align}
The first term in the right-hand side of (\ref{eq:foo45}) is bounded using a Taylor expansion, whereas 
the second term vanishes due to (\ref{eq:L2proj}).  Summing over the elements from $j = 0, \ldots, N$ 
results in
\begin{align}
\label{eq:foo46}
	\int_{0}^L \left(\xi^{n+1} - \xi^n \right) \phi_h  &\leq  C \Delta t^2 \int_0^L |\phi_h|  + \Delta t \sum_{j = 0}^N \left( \mathcal{H}_j(u_h^{n}, \phi_h) - \Delta t  \mathcal{H}_j(u^{n}, \phi_h) \right).
\end{align}
Define
\begin{equation}
\label{eq:defbn1}
\tilde{b}^n(\phi_h) = \Delta t \sum_{j = 0}^N \left( \mathcal{H}_j(u_h^{n}, \phi_h) - \Delta t  \mathcal{H}_j(u^{n}, \phi_h) \right).
\end{equation}
Then equation (\ref{eq:foo46}) becomes
\begin{align}
\label{eq:scheme-scalar_pre}
\int_{0}^L \left(\xi^{n+1} - \xi^n \right) \phi_h  &\leq  C \Delta t^2 \int_0^L |\phi_h|  +\tilde{b}^n(\phi_h),
\end{align}
and Cauchy Schwarz's and Young's inequalities imply
\begin{align}
\label{eq:scheme-scalar}
	\int_{0}^L \left(\xi^{n+1} - \xi^n \right) \phi_h  &\leq  C \Delta t^3 + C \Delta t  \|\phi_h\|^2  +\tilde{b}^n(\phi_h).
\end{align}
We now choose $\phi_h = \xi^n$ to obtain:
\begin{align}
\label{eq:foo47}
\int_{0}^L \left(\xi^{n+1} - \xi^n \right) \xi^n  \leq  C \Delta t^3 + C \Delta t  \|\xi^n\|^2 + \tilde{b}^n(\xi^n).	
\end{align}
It then follows that
\begin{align}
\label{eq:scalar_error_ineq}
	\frac{1}{2} \| \xi^{n+1} \|^2 &- \frac{1}{2} \| \xi^{n} \|^2 \leq \frac{1}{2}\| \xi^{n+1} - \xi^{n} \|^2 +  C \Delta t^3 + C \Delta t  \|\xi^n\|^2 + \tilde{b}^n(\xi^n).
\end{align}
The terms $\Vert \xi^{n+1}-\xi^n\Vert$  and $\tilde{b}^n(\xi^n)$ are bounded by:
\begin{align}
\Vert \xi^{n+1}-\xi^n\Vert^2 \leq C \Delta t^4 + C \Delta t \Vert \xi^n\Vert^2 + C \Delta t\, h^{2k+2}, 
\label{eq:xibn1}\\
\tilde{b}^n(\xi^n) \leq C \Delta t \Vert \xi^n\Vert^2 + C \Delta t \, h^{2k+1}.
\label{eq:xibn2}
\end{align}
Proof of \eqref{eq:xibn1} follows closely the proof of Lemma~\ref{cor:AB1} but is
less technical. We skip it.  Proof of \eqref{eq:xibn2} differs from the proof of Lemma~\ref{cor:AB2} and details are given in
Appendix~\ref{app:proofxibn2}.  
The error inequality simplifies to:
\begin{align*}
 \| \xi^{n+1} \|^2 - \| \xi^{n} \|^2 \leq C \Delta t^3 +  C \Delta t  \|\xi^n\|^2 + C \Delta t\,h^{2k+1}.
\end{align*}
Summing from $n = 0, \ldots, \ell-1$, and using the fact that $\xi^0 = 0$,  one obtains:
\begin{align*}
	\| \xi^{\ell} \|^2 &\leq C \Delta t^2 + C h^{2k+1} + C \Delta t \sum_{n=0}^{\ell-1} \|\xi^n\|^2.
\end{align*}
We now apply Gronwall's inequality:
\begin{align*}
	\| \xi^\ell \|^2 \leq C_4 T \mathrm{e}^T (\Delta t^2 +  h^{2k+1}),
\end{align*}
where $C_4$ is independent of $\ell$.  Employing the CFL condition $\Delta t = O(h^2)$, one has:
\begin{align*}
	\| \xi^\ell \| \leq \left(C_4 T \mathrm{e}^T \right)^{1/2} \left( \Delta t +  h^{k+1/2} \right) 
= \left(C_4 T \mathrm{e}^T \right)^{1/2} \left( h^2 +  h^{k+1/2} \right).
\end{align*} 
Hence the induction is complete if $h$ is small enough so that
\begin{align*}
	C_4 T \mathrm{e}^T h < 1.
\end{align*}
Since $\|\eta^\ell\| \leq C h^{k+1}$ and $\|u^{\ell} - u_h^{\ell}\|  \leq \|\eta^\ell\| + \|\xi^\ell\|$ one obtains:
\begin{align*}
	\| u^\ell - u_h^\ell \| \leq C (\Delta t +  h^{k+1/2}),
\end{align*}
and we conclude the proof.

\section{Numerical results}
\label{sec:numer}

\subsection{Scalar case}

In this section, we use the method of manufactured solutions to numerically verify convergence rates.  Solutions to the inviscid Burger's equation, 
\begin{align}
\label{eq:burgers}
	\frac{\partial u}{\partial t} + \frac{\partial}{\partial x} \left( \frac{1}{2} u^2 \right) = 0,
\end{align}
are approximated using the Adams--Bashforth scheme \eqref{eq:ABscheme}.
We consider the following exact solution to (\ref{eq:burgers}) posed in the interval $[0,1]$:
\begin{align*}
	u(x,t) = \cos(2 \pi x) \sin(t) + \sin(2 \pi x) \cos(t).
\end{align*}
Convergence rates in space, given in Table \ref{table:scalar-rates1}, are calculated for polynomial degrees $k = 1, 2, 3$ by fixing a small timestep $\Delta t = 10^{-4}$  so the temporal error is small compared to the spatial error.  The spatial discretization parameter $h = 1/2^m$ for $m = 1, \ldots, 5$, and we evolve the solution for ten timesteps.  Our results yield a rate of $k+1$ in space, verifying the fact that the convergence estimate in Theorem~\ref{thm:AB} is suboptimal. 

Errors and rates in time are provided in Table \ref{table:scalar-rates2}.  We fix $h = 1/4$, vary $\Delta t = 1/2^m$, $m = 10, \ldots, 13$, and consider high polynomial degrees $k = 8,9$ so the spatial error is smaller than the temporal error.  We evolve the solution to the final time $T = 1$ s.  We recover the expected second order rate in time.
\begin{table}[!htb]
\begin{center}
\scalebox{0.85}{\begin{tabular}{ c | c c | c c | c c }
 & $k=1$ & & $k=2$ & & $k = 3$ &   \\ \hline
$h$ & $L^2$ error & rate & $L^2$ error & rate & $L^2$ error & rate  \\ \hline
5.000$\times 10^{-1}$ & 3.07771$\times 10^{-1}$ & -- & 1.72638$\times 10^{-2}$ & -- & 1.72640$\times 10^{-2}$ & -- \\
2.500$\times 10^{-1}$ & 6.27869$\times 10^{-2}$ & 2.29 & 8.38603$\times 10^{-3}$ & 1.04 & 8.34443$\times 10^{-4}$ & 4.37 \\
1.250$\times 10^{-1}$ & 1.61362$\times 10^{-2}$ & 1.96 & 1.07254$\times 10^{-3}$ & 2.96 & 5.34700$\times 10^{-5}$ & 3.96 \\
6.250$\times 10^{-2}$ & 4.07971$\times 10^{-3}$ & 1.98 & 1.35112$\times 10^{-4}$ & 2.98 & 3.42942$\times 10^{-6}$ & 3.96 \\
3.125$\times 10^{-2}$ & 1.03845$\times 10^{-3}$ & 1.97 & 1.70494$\times 10^{-5}$ & 2.98 & 2.26734$\times 10^{-7}$ & 3.91 \\
\end{tabular}}
\caption{Errors and rates in space for the manufactured solution to Burgers'equation.}
\label{table:scalar-rates1}
\vspace{0.5cm}
\scalebox{0.85}{\begin{tabular}{ c |  c c | c c }
 &  $k=8$ &  & $k = 9$ &   \\ \hline
$\Delta t$ &  $L^2$ error & rate & $L^2$ error & rate  \\ \hline
9.766$\times 10^{-4}$ &  3.01560$\times 10^{-7}$ & -- & 3.04272$\times 10^{-7}$ & -- \\ 
4.883$\times 10^{-4}$ &  7.53310$\times 10^{-8}$ & 2.00 & 7.60427$\times 10^{-8}$ & 2.00 \\ 
2.441$\times 10^{-4}$ &  1.88202$\times 10^{-8}$ & 2.00 & 1.90062$\times 10^{-8}$ & 2.00 \\ 
1.221$\times 10^{-4}$ &  4.87902$\times 10^{-9}$ & 1.94 & 4.74971$\times 10^{-9}$ & 2.00 \\ 
\end{tabular}}
\caption{Errors and rates in time for the manufactured solution to Burgers'equation.}
\label{table:scalar-rates2}
\end{center}
\end{table}

\subsection{System case}

In this section we compute convergence rates for a hyperbolic system that is the motivation for this work: a model which describes one--dimensional blood flow in an elastic vessel:
\begin{align}
\label{eq:AQcons}
&\frac{\partial}{\partial t}\begin{bmatrix}
A \\
Q
\end{bmatrix}
+
\frac{\partial}{\partial x}
\begin{bmatrix}
Q \\
\alpha \frac{Q^2}{A} + \frac{1}{\rho}(A\psi - \Psi)
\end{bmatrix}
=
\begin{bmatrix}
0 \\
-2 \pi \nu \frac{\alpha}{\alpha-1}\frac{Q}{A}
\end{bmatrix}, \\
&p = p_0 + \psi(A; A_0), \quad \Psi = \int_{A_0}^A \psi(\xi;A_0) d \xi.
\end{align}
The variables are vessel cross sectional area $A$ and fluid momentum $Q$.  The parameters are the reference pressure $p_0 = 0$ dynes/$\text{cm}^2$, the reference cross sectional area $A_0 = 1$ $\text{cm}^2$, the non--dimensional Coriolis coefficient $\alpha = 1.1$, the fluid density $\rho = 1.06$ g/$\text{cm}^3$, and the kinematic viscosity $\nu = 3.302 \times 10^{-2}$ $\text{cm}^2$/s.  For these computations we use a typical form for the function relating area to pressure \cite{MGC07}:
\begin{align*}
	\psi = \beta(A^{1/2} - A_0^{1/2}),
\end{align*} 
with $\beta = 1$ dynes/$\text{cm}^3$.  In defining the numerical flux for our computations, we use a version of the local Lax--Friedrichs flux suggested for nonlinear hyperbolic systems in \cite{CS_3}.  With ${\bf U} = [A,Q]^T$ and $\lambda_1({\bf U})$ and $\lambda_2({\bf U})$ the eigenvalues of the Jacobian of the flux function in (\ref{eq:AQcons}), the flux is defined with:
\begin{align*}
	J({\bf U}^-|_{x_j}, {\bf U}^+|_{x_j}) = \max \left( \big|\lambda_1({\bf U}^-|_{x_j})\big|, \big|\lambda_1({\bf U}^+|_{x_j})\big| , \big|\lambda_2({\bf U}^-|_{x_j})\big| , \big|\lambda_2({\bf U}^+|_{x_j})\big|   \right).
\end{align*}

To compute errors and rates, we solve (\ref{eq:AQcons}) in the interval $[0,1]$ with the following exact solution:
\begin{align*}
	A(x,t) = \cos(2 \pi x) \cos(t) + 2,  \quad Q(x,t) = \sin(2 \pi x) \cos(t).
\end{align*}
The discretization for a hyperbolic system follows the same procedure as for a scalar hyperbolic equation.  For these simulations, we employ the second--order Adams--Bashforth scheme  (\ref{eq:ABscheme}) with the local Lax--Friedrichs numerical flux.

Errors and convergence rates in space, provided in Tables \ref{table:rates1} and  \ref{table:rates2} , are determined by fixing a small time step $\Delta t = 2 \times 10^{-5}$ s and taking $h = 1 / 2^m$ for $m = 1, \ldots 5$.  We consider $k = 1, 2, 3$ and evolve the solution for ten time steps.

\begin{table}[!htb]
\begin{center}
\scalebox{0.85}{\begin{tabular}{ c | c c | c c | c c }
 & $k=1$ & & $k=2$ & & $k = 3$ &   \\ \hline
$h$ & $L^2$ error & rate & $L^2$ error & rate & $L^2$ error & rate  \\ \hline
5.000$\times 10^{-1}$ & 8.50463$\times 10^{-2}$ & -- & 8.50463$\times 10^{-2}$ & -- & 2.77383$\times 10^{-3}$ & -- \\
2.500$\times 10^{-1}$ & 6.27702$\times 10^{-2}$ & 0.43 & 8.38200$\times 10^{-3}$ & 3.34 & 8.33345$\times 10^{-4}$ & 1.73 \\
1.250$\times 10^{-1}$ & 1.61152$\times 10^{-2}$ & 1.96 & 1.07125$\times 10^{-3}$ & 2.96 & 5.31039$\times 10^{-5}$ & 3.97 \\
6.250$\times 10^{-2}$ & 4.05695$\times 10^{-3}$ & 1.98 & 1.34722$\times 10^{-4}$ & 2.99 & 3.34118$\times 10^{-6}$ & 3.99 \\
3.125$\times 10^{-2}$ & 1.01713$\times 10^{-3}$ & 1.99 & 1.69031$\times 10^{-5}$ & 2.99 & 2.10357$\times 10^{-7}$ & 3.98 \end{tabular}}
\caption{Errors and rates in space for $A$.}
\label{table:rates1}
\vspace{0.5cm}
\scalebox{0.85}{\begin{tabular}{ c | c c | c c | c c }
 & $k=1$ & & $k=2$ & & $k = 3$ &   \\ \hline
$h$ & $L^2$ error & rate & $L^2$ error & rate & $L^2$ error & rate  \\ \hline
5.000$\times 10^{-1}$ & 3.07761$\times 10^{-1}$ & -- & 1.72654$\times 10^{-2}$ & -- & 1.72638$\times 10^{-2}$ & -- \\
2.500$\times 10^{-1}$ & 6.27688$\times 10^{-2}$ & 2.29 & 8.38233$\times 10^{-3}$ & 1.04 & 8.33176$\times 10^{-4}$ & 4.37 \\
1.250$\times 10^{-1}$ & 1.61145$\times 10^{-2}$ & 1.96 & 1.07130$\times 10^{-3}$ & 2.96 & 5.30850$\times 10^{-5}$ & 3.97 \\
6.250$\times 10^{-2}$ & 4.05679$\times 10^{-3}$ & 1.98 & 1.34717$\times 10^{-4}$ & 2.99 & 3.33998$\times 10^{-6}$ & 3.99 \\
3.125$\times 10^{-2}$ & 1.01736$\times 10^{-3}$ & 1.99 & 1.68933$\times 10^{-5}$ & 2.99 & 2.10567$\times 10^{-7}$ & 3.98 
\end{tabular}}
\caption{Errors and rates in space for $Q$.}
\label{table:rates2}
\end{center}
\end{table}

To calculate the rate in time, we make the error in space small by choosing high order polynomials $k = 8,9$ on a mesh with size $h=1/4$.  By taking $h$ to be constant, we avoid overly refining $\Delta t$ due to the CFL condition.  The time step $\Delta t = 1/2^m$ for $m = 10, \ldots, 13$ and we evolve the solution to the final time $T = 1$ s.  Results are displayed in Tables \ref{table:rates3} and \ref{table:rates4}.

\begin{table}[!htb]
\begin{center}
\scalebox{0.85}{\begin{tabular}{ c |  c c | c c }
 & $k=8$ & & $k = 9$ &   \\ \hline
$\Delta t$ &  $L^2$ error & rate & $L^2$ error & rate  \\ \hline
9.766$\times 10^{-4}$ &  2.90612$\times 10^{-7}$ & -- & 2.98344$\times 10^{-7}$ & -- \\ 
4.883$\times 10^{-4}$ &  7.27141$\times 10^{-8}$ & 1.99 & 7.46399$\times 10^{-8}$ & 1.99 \\ 
2.441$\times 10^{-4}$ &  1.82053$\times 10^{-8}$ & 1.99 & 1.86720$\times 10^{-8}$ & 1.99 \\ 
1.221$\times 10^{-4}$ &  4.59094$\times 10^{-9}$ & 1.98 & 4.67588$\times 10^{-9}$ & 1.99 
\end{tabular}}
\caption{Errors and rates in time for $A$.}
\label{table:rates3}
\vspace{0.5cm}
\scalebox{0.85}{\begin{tabular}{ c |  c c | c c }
 & $k=8$ & & $k = 9$ &   \\ \hline
$\Delta t$ &  $L^2$ error & rate & $L^2$ error & rate  \\ \hline
9.766$\times 10^{-4}$ &  1.88619$\times 10^{-7}$ & -- & 1.91639$\times 10^{-7}$ & -- \\ 
4.883$\times 10^{-4}$ &   4.71556$\times 10^{-8}$ & 1.99 & 4.79006$\times 10^{-8}$ & 2.00 \\ 
2.441$\times 10^{-4}$ &   1.18056$\times 10^{-8}$ & 1.99 & 1.19766$\times 10^{-8}$ & 1.99 \\ 
1.221$\times 10^{-4}$ &   2.99433$\times 10^{-9}$ & 1.97 & 2.99764$\times 10^{-9}$ & 1.99 
\end{tabular}}
\caption{Errors and rates in time for $Q$.}
\label{table:rates4}
\end{center}
\end{table}
 
The computed rates in space and time indicate that results analogous to Theorems \ref{thm:AB} and \ref{thm:FE} can be expected for such numerical discretizations of nonlinear hyperbolic systems.  Numerical analysis for systems will be the subject of future work.

\section{Conclusions}

In this paper we prove a priori error estimates for fully discrete schemes approximating scalar conservation laws, where the spatial discretization is a discontinuous Galerkin method and the temporal discretization is either the second order Adams--Bashforth method or the forward Euler method.  The estimates are valid for polynomial degree greater than or equal to two for the second order method and greater than or equal to one
for the first order method in time. A CFL condition of the form $\Delta t = O(h^2)$ is required.  In future work, we will consider a priori error estimates for numerical methods approximating nonlinear hyperbolic systems like those describing blood flow in an elastic vessel.

\clearpage
\section{Appendix}

\subsection{Proof of bound \eqref{eq:boundF}}
\label{app:boundF}

Using Taylor expansions up to third order, we write
\begin{align}
	f(u_h^n) - f(u^n) &= f'(u^n)(u_h^n - u^n) +\frac12 f''(u^n)(u_h^n-u^n)^2
+\frac16 f'''(\zeta_1^n) (u_h^n-u^n)^3 \nonumber\\
&= f'(u^n)(\chi^n-\eta^n) +\frac12 f''(u^n)(\chi^n-\eta^n)^2 
+\frac16 f'''(\zeta_1^n) (\chi^n-\eta^n)^3 ,  \nonumber\\
&= f'(u^n)\chi^n 
+\frac12 f''(u^n) (\chi^n)^2 
-f'(u^n) \eta^n 
- f''(u^n)  \chi^n \eta^n 
+\frac12 f''(u^n) (\eta^n)^2 
+\frac16 f'''(\zeta_1^n) (\chi^n-\eta^n)^3
\nonumber\\
&= \beta_1+\dots+\beta_6,\nonumber\\
	f(\{u_h^n\}) - f(u^n) &= f'(u^n)(\{u_h^n\} - \{u^n\}) 
+\frac12 f''(u^n)(\{u_h^n\}-u^n)^2\nonumber\\
&=  f'(u^n)(\{\chi^n\} - \{\eta^n\}) + \frac12 f''(u^n)(\{\chi^n\}-\{\eta^n\})^2
+\frac16 f'''(\zeta_2^n) (\{\chi^n\}-\{\eta^n\})^3,\nonumber\\
&= f'(u^n)\{\chi^n\} 
+\frac12 f''(u^n) (\{\chi^n\})^2 
-f'(u^n) \{\eta^n\} 
- f''(u^n) \{\chi^n\} \{\eta^n\} 
\nonumber\\
&
+\frac12 f''(u^n) (\{\eta^n\})^2 
+\frac16 f'''(\zeta_2^n) (\{\chi^n\}-\{\eta^n\})^3
\nonumber\\
&= \gamma_1+\dots+\gamma_6,
\end{align}
where $\zeta^n_1$ and $\zeta_2^n$ are some points between $u_h^n$ and $u^n$, and $\{u_h^n\}$ and $u^n$ respectively.
We substitute these expansions in the terms $\mathcal{F}(n,\chi^n)$ and write:
\begin{align}
        \mathcal{F}(n,\chi^n) &= X_1 + \ldots + X_6,
\end{align}
with
\begin{align}
        X_i = \Delta t \sum_{j = 0}^N  \int_{I_j}\beta_i\frac{d \chi^n}{dx} 
- \Delta t \sum_{j = 1}^N \gamma_i |_{x_j} [\chi^n]|_{x_j}, \quad 1\leq i\leq 6.
\end{align}
We integrate by parts the first term in the definition of $X_1$ and use the fact that $f'$ vanishes at the endpoints of the domain, namely
at $x_0$ and $x_{N+1}$. The term $X_1$ then simplifies to
\[
X_1 =  -\frac12 \Delta t \sum_{j=0}^N \int_{I_j} (\frac{\partial}{\partial x} f'(u^n)) (\chi^n)^2
\leq C \Delta t \Vert \chi^n \Vert^2.
\]
Using the assumption $\Vert \chi^n\Vert \leq h^{3/2}$ and trace inequalities, we have
\begin{align}
X_2  &= \frac12 \Delta t \sum_{j=0}^N \int_{I_j} f''(u^n) (\chi^n)^2 \frac{d \chi^n}{dx}
-\frac12 \Delta t \sum_{j=1}^N f''(u^n)|_{x_j}  (\{\chi^n\})^2|_{x_j} [\chi^n]|_{x_j}\nonumber\\
&\leq C \Delta t \Vert \chi^n\Vert_{\infty} h^{-1} \Vert \chi^n\Vert^2
\leq C \Delta t  \Vert \chi^n\Vert^2.
\end{align}
To bound the term $X_3$ we define the following piecewise
constant function $u_c^n$ elementwise as:
\begin{align}
	u_c^n|_{I_j}(x) = u^{n}|_{x_j}, \quad \forall x \in I_j, \quad \forall 0\leq j \leq N.
\end{align}
We  note that 
\begin{align}
\label{eq:foo21}
	\|f'(u^n) - f'(u_c^n)\|_\infty \leq Ch. 
\end{align}
We then rewrite the term $X_3$
\begin{align}
	X_3 & =  -\Delta t \sum_{j = 0}^N \int_{I_j} f'(u^n)\eta^n  \frac{d \chi^n}{dx} 
+ \Delta t \sum_{j = 1}^N  f'(u^n)\{\eta^n\} |_{x_j} [\chi^n]|_{x_j}\nonumber\\
& = -\Delta t \sum_{j = 0}^N \int_{I_j} (f'(u^n)-f'(u_c^n))\eta^n  \frac{d \chi^n}{dx} 
-\Delta t \sum_{j = 0}^N f'(u_c^n) \int_{I_j} \eta^n  \frac{d \chi^n}{dx} \nonumber\\
&+ \Delta t \sum_{j = 1}^N  (f'(u^n)-f'(\{u_h^n\}))|_{x_j} \{\eta^n\} |_{x_j} [\chi^n]|_{x_j}
+ \Delta t \sum_{j = 1}^N  f'(\{u_h^n\})|_{x_j} \{\eta^n\} |_{x_j} [\chi^n]|_{x_j}.
\label{eq:interX3}
\end{align}
The second term above vanishes because of \eqref{eq:L2proj}.
The first term is bounded using approximation properties and \eqref{eq:foo21}.
\[
\Delta t \sum_{j = 0}^N \int_{I_j} (f'(u^n)-f'(u_c^n))\eta^n  \frac{d \chi^n}{dx}\leq C \Delta t\, h^{2k+2} +
C\Delta t \Vert \chi^n\Vert^2.
\]
Using a Taylor expansion, for some $\zeta_3^n$ we have
\[
f'(u^n)-f'(\{u_h^n\}) = f''(\zeta_3^n) \{u^n-u_h^n\} \leq C (\Vert \chi^n\Vert_\infty + \Vert \eta^n\Vert_\infty).
%HOW DO WE KNOW WE HAVE LINFTY BOUND FOR ETA?
\]
Using the assumption $\Vert \chi^n\Vert \leq h^{3/2}$ we then have
\[
\Delta t \sum_{j = 1}^N  (f'(u^n)-f'(\{u_h^n\}))\{\eta^n\} |_{x_j} [\chi^n]|_{x_j}
\leq C \Delta t \Vert \chi^n\Vert^2 + C \Delta t \, h^{2k+2}.
\]
For the last term in (\ref{eq:interX3}) we employ (\ref{eq:alpha1}) to obtain:
\begin{align}
\Delta t \sum_{j = 1}^N  &f'(\{u_h^n\})|_{x_j} \{\eta^n\} |_{x_j} [\chi^n]|_{x_j} 
\leq C\Delta t \sum_{j = 1}^N \left(\alpha(u_h^n)|_{x_j}  + C |[u_h^n]|_{x_j} | \right) \Big|\{\eta^n\} |_{x_j} \Big| \Big| [\chi^n]|_{x_j} \Big| \nonumber \\
	&= C\Delta t \sum_{j = 1}^N\alpha(u_h^n)|_{x_j}  \Big|\{\eta^n\} |_{x_j} \Big| \Big| [\chi^n]|_{x_j} \Big| 
+ C\Delta t \sum_{j = 1}^N |[u_h^n]|  \Big|\{\eta^n\} |_{x_j} \Big| \Big| [\chi^n]|_{x_j} \Big| \label{eq:foo23}. 
\end{align}
Using Cauchy-Schwarz's and Young's inequalities, approximation results and the assumption $\Vert \chi^n\Vert \leq h^{3/2}$, we obtain
\[
\Delta t \sum_{j = 1}^N  f'(\{u_h^n\})|_{x_j} \{\eta^n\} |_{x_j} [\chi^n]|_{x_j} 
\leq \varepsilon \Delta t \sum_{j=1}^N \alpha(u_h^n)|_{x_j}  [\chi^n]|_{x_j}^2
+ C\varepsilon^{-1} \Delta t \, h^{2k+1}
+ C \Delta t \Vert \chi^n\Vert^2.
\]
In summary we have
\[
X_3 \leq \varepsilon \Delta t \sum_{j=1}^N \alpha(u_h^n)|_{x_j}  [\chi^n]|_{x_j}^2
+ C\varepsilon^{-1} \Delta t \, h^{2k+1}
+ C \Delta t \Vert \chi^n\Vert^2 + C \Delta t\, h^{2k+1}.
\]
The bounds for $X_4$, $X_5$, and $X_6$ are standard applications of Cauchy Schwarz's inequality, Young's inequality, the induction hypothesis, assumption (\ref{eq:fbounded}), and inequalities (\ref{eq:inv2}), (\ref{eq:inv3}), and (\ref{eq:est1})--(\ref{eq:est3}):
\begin{align}
X_4 & = -\Delta t \sum_{j=0}^N \int_{I_j} f''(u_n) \chi^n \eta^n \frac{d \chi^n}{dx}
+\Delta t \sum_{j=1}^N f''(u_n)|_{x_j}  \{\chi^n\}|_{x_j} \{\eta^n\}|_{x_j} [\chi^n]|_{x_j}
\nonumber\\
	 &\leq C \Delta t \|\chi^n\|^2 + \Delta t\,h^{2k+1}, \\
X_5 & = \frac12 \Delta t \sum_{j=0}^N \int_{I_j} f''(u_n)|_{x_j}  (\eta^n)^2 \frac{d \chi^n}{dx}
-\frac12 \Delta t \sum_{j=1}^N f''(u_n) \{\eta^n\}^2|_{x_j}  [\chi^n]|_{x_j}
\nonumber\\
	&\leq C \Delta t\,h^{2k+2} + C \Delta t \|\chi^n\|^2, \\
X_6 &= \frac16 \Delta t \sum_{j=0}^N \int_{I_j} f'''(\zeta_1^n) (\chi^n-\eta^n)^3 \frac{d \chi^n}{dx}
-\frac16 \Delta t \sum_{j=1}^N f'''(\zeta_2^n) (\{\chi^n\}-\{\eta^n\})^3|_{x_j} [\chi^n]|_{x_j}
\nonumber\\
	 &\leq C \Delta t \|\chi^n\|^2 + C \Delta t\,h^{2k+1}.
\end{align}
We can then conclude by combining all the bounds above. 

\subsection{Proof of bound \eqref{eq:xibn2}}
\label{app:proofxibn2}

We rewrite, using the definition of $\alpha$
\[
\tilde{b}^n(\xi^n) = \theta_1 + \theta_2 + \theta_3,
\]
with
\[
\theta_1 = 
\Delta t \sum_{j=0}^N \int_{I_j} (f(\tilde{u}_h^n)-f(u^n))\frac{d\xi^n}{dx}
- \Delta t \sum_{j=1}^N (f(\{\tilde{u}_h^n\})-f(u^n))|_{x_j} [\xi^n]|_{x_j},
\]
\[
\theta_2 =  \Delta t \sum_{j=1}^N \int_{I_j} (s(\tilde{u}_h^n)-s(u^n))\xi^n,
\]
\[
\theta_3 = -\Delta t \sum_{j=1}^N \alpha(\tilde{u}_h^n)|_{x_j} [\tilde{u}_h^n]|_{x_j} [\xi^n]|_{x_j}.
\]
We note that the bound for $\theta_1$
follows the argument of the proof of \eqref{eq:boundF}, where we substitute $\chi^n$ by $\xi^n$.
As in the previous section,  we use Taylor expansions up to third order and write the term $\tilde{b}^n(\xi^n)$ as a sum
of six terms, $X_i$, $1\leq i\leq 6$.  Bounds for $X_i$ are obtained in a similar fashion, except for the term $X_2$
which is bounded differently because the 
the induction hypothesis for the forward Euler scheme is weaker than the hypothesis for the Adams--Bashforth scheme.
We have
\begin{align}
        X_2 &=  \Delta t \frac{1}{2} \sum_{j = 0}^N  \int_{I_j}f''(u^n)(\xi^n)^2\frac{d \xi^n}{dx} - \Delta t \frac{1}{2} \sum_{j = 1}^N f''(u^n)(\{\xi^n\})^2 |_{x_j} [\xi^n]|_{x_j}.
\end{align}
We rewrite the first term above. Integrating the first term by parts gives and using the assumption that $f'''$ vanishes at the endpoints of the interval gives:
\begin{align}
        \Delta t \frac{1}{2} \sum_{j = 0}^N  \int_{I_j}f''(u^n)(\xi^n)^2\frac{d \xi^n}{dx} &= \Delta t \frac{1}{6} \sum_{j = 0}^N  \int_{I_j}f''(u^n)\frac{d (\xi^n)^3}{dx} \nonumber\\
        \label{eq:foo20}
        &= \Delta t \frac{1}{6} \sum_{j = 1}^N f''(u^n)|_{x_j} [(\xi^n)^3]|_{x_j} -\Delta t \frac{1}{6} \sum_{j = 0}^N \int_{I_j}\frac{\partial f''(u^n)}{\partial x} (\xi^n)^3.
\end{align}
Now, we use the identity $[\xi^3] = 2\{\xi\}^2[\xi] + \{ \xi^2\}[\xi]$ to rewrite the first term in the right-hand side of \eqref{eq:foo20}:
\begin{align}
        X_2 &= \Delta t \frac{1}{6} \sum_{j = 1}^N f''(u^n)|_{x_j} \left( \{(\xi^n)^2\} - \{\xi^n\}^2\right) [\xi^n]|_{x_j} - \Delta t \frac{1}{6} \sum_{j = 0}^N \int_{I_j}\frac{\partial f''(u^n)}{\partial x} (\xi^n)^3.
\end{align}
Employing the identity $\{\xi^2\} - \{\xi\}^2  = \frac{1}{4}[\xi]^2$ for the first term and inductive hypothesis $\|\xi^n\|_\infty \leq h^{1/2}$ on the second term gives:
\begin{align}
        X_2 &\leq \Delta t \frac{1}{24} \sum_{j = 1}^N f''(u^n)|_{x_j}  [\xi^n]^3|_{x_j} + C \Delta t \|\xi\|_\infty \|\xi\|^2 \leq \Delta t \frac{1}{24} \sum_{j = 1}^N f''(u^n)|_{x_j}  [\xi^n]^3|_{x_j} + C \Delta t \|\xi\|^2.
        \label{eq:foo85}
\end{align}
The first term in \eqref{eq:foo85} is broken into two parts:
\begin{align}
\label{eq:foo86}
        \Delta t \frac{1}{24} \sum_{j = 1}^N f''(u^n)|_{x_j}  [\xi^n]^3|_{x_j}  &= \Delta t \frac{1}{24}  \sum_{j = 1}^N \left(f''(u^n)|_{x_j} - f''(\{\tilde{u}_h^n\})|_{x_j} \right)  [\xi^n]^3|_{x_j} +   \Delta t \frac{1}{24} \sum_{j = 1}^N f''(\{\tilde{u}_h^n\})|_{x_j} [\xi^n]^3|_{x_j}.
\end{align}
We use for the first term in (\ref{eq:foo86}) a Taylor expansion $f''(u^n) - f''(\{\tilde{u}_h^n\}) = f'''(\zeta^n) \{ \eta^n-\xi^n\}$ 
with the inductive hypothesis to obtain the following bound:
\begin{align}
        \Delta t \frac{1}{24} \sum_{j = 1}^N (f''(u^n)|_{x_j}-f''(u^n_h)|_{x_j})  [\xi^n]^3|_{x_j}
%&\leq C \Delta t \|e^n\|_\infty \|\xi^n\|_\infty \sum_{j = 1}^N [\xi^n]^2|_{x_j} + \Delta t \frac{1}{24} \sum_{j = 1}^N f''(\{u_h^n\})|_{x_j} [\xi^n]^3|_{x_j} \nonumber \\
        \label{eq:foo87}
        \leq C\Delta t \|\xi^n\|^2.
\end{align}
For the last term in (\ref{eq:foo86}), since $[u^n] = 0$, we rewrite it using the identity $[\xi^n] = [\eta^n] + [\tilde{u}_h^n]$:
\begin{align}
\label{eq:foo88}
        \Delta t \frac{1}{24} \sum_{j = 1}^N f''(\{\tilde{u}_h^n\})|_{x_j} [\xi^n]^3|_{x_j} & 
= \Delta t \frac{1}{24} \sum_{j = 1}^N f''(\{\tilde{u}_h^n\})|_{x_j} [\eta^n] |_{x_j}[\xi^n]^2|_{x_j} 
+ \Delta t \frac{1}{24} \sum_{j = 1}^N f''(\{\tilde{u}_h^n\})|_{x_j} [\tilde{u}_h^n]|_{x_j} [\xi^n]^2|_{x_j}.
        \end{align}
The first term in (\ref{eq:foo88}) can be estimated with trace inequalities and approximation results. The second term in (\ref{eq:foo88}) is bounded using inequality (\ref{eq:alpha2}) and the induction hypothesis:
\begin{align}
        \Delta t \frac{1}{24} \sum_{j = 1}^N f''(\{\tilde{u}_h^n\})|_{x_j} [\xi^n]^3|_{x_j} &\leq C \Delta t h^{-1} \|\eta^n\|_\infty \|\xi^n\|^2 
+ \Delta t \frac{1}{3} \sum_{j = 1}^N \left(\alpha(\tilde{u}_h^n)|_{x_j} + C |[\tilde{u}_h^n]|^2|_{x_j} \right) [\xi^n]^2|_{x_j} \nonumber\\
        &\leq C \Delta t \|\xi^n\|^2 + \Delta t \frac{1}{3} \sum_{j = 1}^N\alpha(\tilde{u}_h^n)|_{x_j}  [\xi^n]^2|_{x_j} + C \Delta t\,h^{-1}\|\tilde{u}_h^n\|_\infty^2 \|\xi^n\|^2 \nonumber\\
        &\leq C \Delta t \|\xi^n\|^2 + \Delta t \frac{1}{3} \sum_{j = 1}^N\alpha(\tilde{u}_h^n)|_{x_j}  [\xi^n]^2|_{x_j}.
\end{align}
Combining all the estimates gives:
\begin{align}
X_2 \leq \Delta t \frac{1}{3} \sum_{j = 1}^N\alpha(\tilde{u}_h^n)|_{x_j}  [\xi^n]^2|_{x_j} + C\Delta t\|\xi^n\|^2.
\end{align}
This bound is added to the bounds for the other terms $X_i$'s to obtain:
\[
\theta_1 \leq C \Delta t\, \Vert\xi^n\Vert^2 
+ (\frac13 + \varepsilon) \sum_{j = 1}^N\alpha(\tilde{u}_h^n)|_{x_j} [\xi^n]^2|_{x_j} 
+ C \Delta t (1+\varepsilon^{-1}) h^{2k+1}.
\]
The term $\theta_2$ is bounded using Lischitz continuity of $s$:
\[
\theta_2 \leq C \Delta t \Vert \xi^n\Vert^2 + C \Delta t \, h^{2k+2}.
\]
The term $\theta_3$ is rewritten as
\[
\theta_3 = -\Delta t \sum_{j=1}^N \alpha(\tilde{u}_h^n)|_{x_j} [\xi^n]^2|_{x_j}
+ \Delta t \sum_{j=1}^N \alpha(\tilde{u}_h^n)|_{x_j} [\eta^n]|_{x_j} [\xi^n]|_{x_j}.
\]
Using Young's inequality and approximation results we obtain
\[
\theta_3 \leq (-1+\varepsilon) \Delta t \sum_{j=1}^N \alpha(\tilde{u}_h^n)|_{x_j} [\xi^n]^2|_{x_j}
+ C \Delta t\, h^{2k+1}, \quad \forall \varepsilon>0.
\]
This means that by choosing $\varepsilon = 1/3$ in the above, we conclude
\[
\tilde{b}^n(\xi^n) \leq C \Delta t \Vert \xi^n\Vert^2 + C \Delta t \, h^{2k+1}.
\]

\bibliographystyle{plain}
\bibliography{rkdg_refs}

\begin{thebibliography}{10}

\bibitem{Alastruey08}
J.~Alastruey, S.M. Moore, K.H. Parker, T.~David, J.~Peir{\'o}, and S.J.
  Sherwin.
\newblock Reduced modelling of blood flow in the cerebral circulation: coupling
  1-{D}, 0-{D} and cerebral auto-regulation models.
\newblock {\em Internat. J. Numer. Methods Fluids}, 56(8):1061, 2008.

\bibitem{Boi15}
E.~Boileau, P.~Nithiarasu, P.J. Blanco, L.O. M{\"u}ller, F.E. Fossan, L.R.
  Hellevik, W.P. Donders, W.~Huberts, M.~Willemet, and J.~Alastruey.
\newblock A benchmark study of numerical schemes for one-dimensional arterial
  blood flow modelling.
\newblock {\em Internat. J. Numer. Methods Biomed. Eng.}, 2015.

\bibitem{Bollache2014}
E.~Bollache, N.~Kachenoura, A.~Redheuil, F.~Frouin, E.~Mousseaux, P.~Recho, and
  D.~Lucor.
\newblock Descending aorta subject-specific one-dimensional model validated
  against in vivo data.
\newblock {\em Journal of Biomechanics}, 47(2):424--431, 2014.

\bibitem{CK03}
S.~{\v{C}}ani{\'c} and E.H. Kim.
\newblock Mathematical analysis of the quasilinear effects in a hyperbolic
  model blood flow through compliant axi-symmetric vessels.
\newblock {\em Math. Methods Appl. Sci.}, 26(14):1161--1186, 2003.

\bibitem{Cascaval2016}
R.C. Cascaval, C.~D'Apice, M.P. D'Arienzo, and R.~Manzo.
\newblock Boundary control for an arterial system.
\newblock {\em Journal of Fluid Flow}, 3, 2016.

\bibitem{CS_4}
B.~Cockburn, S.~Hou, and C.W. Shu.
\newblock The {R}unge-kutta local projection discontinuous {G}alerkin finite
  element method for conservation laws. {IV}. the multidimensional case.
\newblock {\em Math. Comp.}, 54(190):545--581, 1990.

\bibitem{CS_3}
B.~Cockburn, S.Y. Lin, and C.W. Shu.
\newblock {TVB} {R}unge-{K}utta local projection discontinuous {G}alerkin
  finite element method for conservation laws {III}: one-dimensional systems.
\newblock {\em J. Comput. Phys.}, 84(1):90--113, 1989.

\bibitem{CS_2}
B.~Cockburn and C.W. Shu.
\newblock {TVB} {R}unge-{K}utta local projection discontinuous {G}alerkin
  finite element method for conservation laws. {II}. general framework.
\newblock {\em Math. Comp.}, 52(186):411--435, 1989.

\bibitem{CS_1}
B.~Cockburn and C.W. Shu.
\newblock The {R}unge-{K}utta local projection {$P^1$} discontinuous-{G}alerkin
  finite element method for scalar conservation laws.
\newblock {\em RAIRO-Mod\'elisation {M}ath\'ematique et {A}nalyse
  {N}um\'erique}, 25(3):337--361, 1991.

\bibitem{CS_5}
B.~Cockburn and C.W. Shu.
\newblock The {R}unge-{K}utta discontinuous {G}alerkin method for conservation
  laws {V}: multidimensional systems.
\newblock {\em J. Comput. Phys.}, 141(2):199--224, 1998.

\bibitem{Dafermos2010}
C.M. Dafermos.
\newblock {\em Hyperbolic conservation laws in continuum physics, volume 325 of
  Grundlehren der Mathematischen Wissenschaften [Fundamental Principles of
  Mathematical Sciences]}.
\newblock Springer-Verlag, Berlin, 2010.

\bibitem{Deriaz2012}
E.~Deriaz.
\newblock Stability conditions for the numerical solution of
  convection-dominated problems with skew-symmetric discretizations.
\newblock {\em SIAM J. Numer. Anal.}, 50(3):1058--1085, 2012.

\bibitem{Di2011}
D.A. Di~Pietro and A.~Ern.
\newblock {\em Mathematical aspects of discontinuous Galerkin methods},
  volume~69.
\newblock Springer Science \& Business Media, 2011.

\bibitem{Dumas2017}
L.~Dumas, T.~El~Bouti, and D.~Lucor.
\newblock A robust and subject-specific hemodynamic model of the lower limb
  based on noninvasive arterial measurements.
\newblock {\em Journal of {B}iomechanical {E}ngineering}, 139(1):011002, 2017.

\bibitem{LSZ15}
J.~Luo, C.W. Shu, and Q.~Zhang.
\newblock A priori error estimates to smooth solutions of the third order
  {R}unge--{K}utta discontinuous {G}alerkin method for symmetrizable systems of
  conservation laws.
\newblock {\em ESAIM: Math. Model. Numer. Anal.}, 49(4):991--1018, 2015.

\bibitem{Matthys2007}
K.S. Matthys, J.~Alastruey, J.~Peir{\'o}, A.W. Khir, P.~Segers, P.R. Verdonck,
  K.H. Parker, and S.J. Sherwin.
\newblock Pulse wave propagation in a model human arterial network: assessment
  of 1-{D} numerical simulations against in vitro measurements.
\newblock {\em Journal of Biomechanics}, 40(15):3476--3486, 2007.

\bibitem{MGC07}
A.~Mikelic, G.~Guidoboni, and S.~{\v{C}}ani{\'c}.
\newblock Fluid-structure interaction in a pre-stressed tube with thick elastic
  walls {I}: the stationary {S}tokes\ problem.
\newblock {\em Netw. Heterog. Media}, 2(3):397, 2007.

\bibitem{PCRR16}
C.~Puelz, B.~Rivi{\`e}re, S.~\v{C}ani\'c, and C.G. Rusin.
\newblock Comparison of reduced blood flow models using {R}unge--{K}utta
  discontinuous {G}alerkin methods.
\newblock {\em Appl. Numer. Math.}, 115:114--141, 2017.

\bibitem{SFPF03}
S.J. Sherwin, L.~Formaggia, J.~Peiro, and V.~Franke.
\newblock Computational modelling of \textsc{1D} blood flow with variable
  mechanical properties and its application to the simulation of wave
  propagation in the human arterial system.
\newblock {\em Internat. J. Numer. Methods Fluids}, 43(6-7):673--700, 2003.

\bibitem{SFPP03}
S.J. Sherwin, V.~Franke, J.~Peir{\'o}, and K.~Parker.
\newblock One-dimensional modelling of a vascular network in space-time
  variables.
\newblock {\em J. {E}ngrg. {M}ath.}, 47(3-4):217--250, 2003.

\bibitem{Wang2016}
H.~Wang, C.W. Shu, and Q.~Zhang.
\newblock Stability analysis and error estimates of local discontinuous
  {G}alerkin methods with implicit--explicit time-marching for nonlinear
  convection--diffusion problems.
\newblock {\em Appl. Math. Comput.}, 272:237--258, 2016.

\bibitem{WFL14}
X.~Wang, J.M. Fullana, and P.Y. Lagr{\'e}e.
\newblock Verification and comparison of four numerical schemes for a 1{D}
  viscoelastic blood flow model.
\newblock {\em Computer Methods in Biomechanics and Biomedical Engineering},
  18(15):1704--1725, 2015.

\bibitem{Zakerzadeh16}
M.~Zakerzadeh and G.~May.
\newblock On the convergence of a shock capturing discontinuous galerkin method
  for nonlinear hyperbolic systems of conservation laws.
\newblock {\em SIAM J. Numer. Anal.}, 54(2):874--898, 2016.

\bibitem{ZS04}
Q.~Zhang and C.W. Shu.
\newblock Error estimates to smooth solutions of {R}unge-{K}utta discontinuous
  {G}alerkin methods for scalar conservation laws.
\newblock {\em SIAM J. Numer. Anal.}, 42(2):641--666, 2004.

\bibitem{ZS06}
Q.~Zhang and C.W. Shu.
\newblock Error estimates to smooth solutions of {R}unge-{K}utta discontinuous
  {G}alerkin method for symmetrizable systems of conservation laws.
\newblock {\em SIAM J. Numer. Anal.}, 44(4):1703--1720, 2006.

\bibitem{ZS10}
Q.~Zhang and C.W. Shu.
\newblock Stability analysis and a priori error estimates of the third order
  explicit {R}unge-{K}utta discontinuous {G}alerkin method for scalar
  conservation laws.
\newblock {\em SIAM J. Numer. Anal.}, 48(3):1038--1063, 2010.

\end{thebibliography}
\end{document}